\documentclass[fleqn, a4paper,11pt]{article}

\usepackage{amsmath,txfonts,bm}
\usepackage{amsfonts}
\usepackage{graphicx}
\usepackage{theorem}
\theorembodyfont{\upshape}
\theoremstyle{break}
\newtheorem{thm}{Theorem}

\newtheorem{lem}[thm]{Lemma}
\newtheorem{rmk}[thm]{Remark}

\newtheorem{dfn}[thm]{Definition}

\begin{document}

    \noindent 
    \LARGE{A unified understanding of the two formulae for the traces of the inverse powers of a positive definite symmetric tridiagonal matrix}  \\

    \noindent 
    \large{Takumi Yamashita
           \footnote{28-20 Kojogaoka, Otsu, Shiga 520-0821 Japan, Tel.: +81-77-522-7447, Fax:  +81-77-522-7447, e-mail: t-yamashita@kke.biglobe.ne.jp}
          }

\section*{Abstract}

    For an upper bidiagonal matrix $B$ where all the diagonal and the upper subdiagonal entries are positive,
    two subtraction-free formulae for computation of the traces $J_{M} ( B ) = \textrm{Tr} ( ( B^{\top} B )^{- M} ) = \textrm{Tr} ( ( B B^{\top} )^{- M} )$ $( M = 1, 2, \dots )$
    have been presented in the two preceding works.
    A few lower bounds of the minimal singular value of $B$ are obtained from these traces.
    In this paper, we clarify  some properties of these formulae and present a new subtraction-free formula for the traces $J_{M} ( B )$.
    An interpretation of       some quantities in one of the preceding works in terms of matrix theory is given.
    Some relationships between some quantities in the preceding works are also given.
    From these relationships, the new subtraction-free formula for the traces $J_{M} ( B )$ is obtained.

\section{Introduction}

    Singular values of a matrix are theoretically important and have wide applications.
    Their exact values are, needless to say, have a lot of examples of applications.
    In addition to such examples, estimation of a lower bound of the minimal singular value also has examples of applications.
    An example is estimation of an upper bound of condition number of a matrix.
    In a system of simultaneous linear equations or an eigenvalue problem,
    condition number of the coefficient matrix is used as an indicator for reliability of numerical computation.
    Condition number $\kappa ( A )$ of a matrix $A$ is defined by norms of the matrices $A$ and $A^{- 1}$ as $\kappa ( A ) = \| A \| \cdot \| A^{- 1} \|$.
    When we choose the 2-norm for definition of condition number,
    this value is expressed by     the maximal singular value $\sigma _{\max} ( A )$
                               and the minimal singular value $\sigma _{\min} ( A )$ of $A$ as $\kappa ( A ) = \sigma _{\max} ( A ) / \sigma _{\min} ( A )$.
    The exact values of the maximal and the minimal singular values of $A$ are usually unknown before solving the problem.
    Then, an upper bound of condition number is used as an alternative in practice.
    An upper bound of $\sigma _{\max} ( A )$ is given as $\sigma _{\max} ( A ) \leq \sqrt{\| A \| _{1} \cdot \| A \| _{\infty}}$,
    where $\| A \| _{1}$ and $\| A \| _{\infty}$ are the 1-norm and the infinity norm of $A$, respectively \cite{Johnson89}.
    Lower bounds of the minimal singular value of a matrix have been proposed by several authors.
    See \cite{HP92, Johnson89, JS98, Li01, Rojo99, YG97}, for example.
    As another example of applications, a lower bound of the minimal singular value of a bidiagonal matrix can be applied to singular value computing.
    In standard procedure for singular value computing,
    the input matrix is firstly transformed into a bidiagonal matrix by orthogonal transformations without changing its singular values.
    Then, we obtain the singular values by applying some algorithm to this bidiagonal matrix.
    For example, by     the dqds (the differential quotient-difference algorithm with shifts) algorithm \cite{FP94},
                        the orthogonal qd-algorithm \cite{vonMatt97},
                    or the mdLVs algorithm (modified discrete Lotka-Volterra algorithm with shift) \cite{IN06},
    we can obtain singular values of a bidiagonal matrix.
    These algorithms are iterative and equip a technique called the shift of origin.
    This technique accelerates convergence of iteration.
    In this technique, we are required to give a quantity called a shift at each iteration.
    A too large shift causes failure of the algorithms.
    A shift can be determined from a lower bound of the minimal singular value of the bidiagonal matrix.
    Kimura {\it et al}. \cite{KYN11} presented a sequence of lower bounds of the minimal singular value $\sigma _{\min} ( B )$ of an upper bidiagonal matrix $B$,
    where all the diagonal and the upper subdiagonal entries are positive.
    Let these lower bounds be denoted by $\theta _{M} ( B )$ for $M = 1, 2, \dots$.
    They are given as
    \begin{equation*}
        \theta _{M} ( B ) = \left( J_{M} ( B ) \right) ^{- \frac{1}{2 M}},
    \end{equation*}
    where $J_{M} ( B ) = \textrm{Tr} \left( ( B^{\top} B )^{- M } \right)$.
    Note that it holds $\textrm{Tr} \left( ( B^{\top} B )^{- M } \right) = \textrm{Tr} \left( ( B B^{\top} )^{- M } \right)$.
    These lower bounds increase monotonically and converges to the minimal singular value $\sigma _{\min} ( B )$ as $M \rightarrow \infty$ \cite[Theorem 3.1]{KYN11}, namely,
    \begin{gather*}
                            \theta _{1} ( B ) < \theta _{2} ( B ) < \cdots < \sigma _{\min} ( B ),  \\
        \lim_{M \to \infty} \theta _{M} ( B ) =                              \sigma _{\min} ( B ).
    \end{gather*}
    Kimura {\it et al}. \cite{KYN11} presented also a formula for computation of the trace $J_{M} ( B )$ for a fixed positive integer $M$.
    This formula gives the diagonal entries of $( B^{\top} B )^{- M}$ and $( B B^{\top} )^{- M}$ in a form of recurrence relation.
    In the case of $M \geq 2$, the recurrence relation includes subtraction in it.
    Subtraction may occur cancellation error in numerical computation.
    To exclude this possibility,
    Yamashita, Kimura and Nakamura \cite{YKN12} derived a ``subtraction-free'' recurrence formula,
    which consists of only summation, multiplication and division among positive quantities,
    for the diagonal entries of $( B^{\top} B )^{- M}$ and $( B B^{\top} )^{- M}$ $( M = 1, 2, \dots )$ starting from the formula in \cite{KYN11}.
    On the other hand, Yamashita, Kimura and Yamamoto \cite{YKY14} derived another subtraction-free formula for the traces $J_{M} ( B )$ $( M = 1, 2, \dots )$
    with an idea which is quite different from that in \cite{YKN12}.
    This derivation does not aim to obtain the diagonal entries of $( B^{\top} B )^{- M}$ or $( B B^{\top} )^{- M}$ $( M = 1, 2, \dots )$.
    These formulae in \cite{KYN11}, \cite{YKN12} and \cite{YKY14} are briefly reviewed in the next section.

    In this paper, we consider the two formulae in the preceding works \cite{YKN12} and \cite{YKY14} and achieve a unified understanding of them.
    Firstly,  we show that some quantities used in the formula in \cite{YKN12} can be interpreted in terms of matrix theory.
    Secondly, we reveal relationships between some quantities in the two formulae.
    Lastly, a new subtraction-free formula for the traces $J_{M} ( B )$ $( M = 1, 2, \dots )$ is obtained from these relationships.

    This paper is organized as follows.
    In Section \ref{SecPrecWorks}, we briefly review the formulae in the preceding works \cite{KYN11}, \cite{YKN12} and \cite{YKY14}.
    In Section \ref{NewFormula}, we achieve the unified understanding of the formulae in \cite{YKN12} and \cite{YKY14} and present the new formula.
    Section \ref{ConcRmks} is devoted for concluding remarks.

\section{The formulae for the traces in the preceding works}  \label{SecPrecWorks}

    In this section, we briefly review the formulae for the traces $J_{M} ( B )$ $( M = 1, 2, \dots )$ in the preceding works \cite{KYN11}, \cite{YKN12} and \cite{YKY14}.
    The review of \cite{KYN11} helps to explain the formula in \cite{YKN12}.
    Before review, notation of the $N \times N$ upper bidiagonal matrix $B = ( B_{i, j} )$, where all the diagonal and the upper subdiagonal entries are positive, is necessary.
    Let the matrix $B$ be denoted by
    \begin{equation}  \label{UpperBidiagonal}
        B = \left( \begin{array}{ccccc}
                       \sqrt{q_{1}}  &  \sqrt{e_{1}}  &                 &                    &                    \\
                                     &  \sqrt{q_{2}}  &  \sqrt{e_{2}}   &                    &                    \\
                                     &                &  \ddots         &  \ddots            &                    \\
                                     &                &                 &  \sqrt{q_{N - 1}}  &  \sqrt{e_{N - 1}}  \\
                                     &                &                 &                    &  \sqrt{q_{N    }}
                   \end{array}                                                                                       \right) ,
    \end{equation}
    where $q_{i} > 0$ for $i = 1, \dots , N$ and $e_{i} > 0$ for $i = 1, \dots , N - 1$.
    When we consider singular values of an upper bidiagonal matrix,
    we can impose the condition of positivity of the diagonal and the upper subdiagonal entries without loss of generality \cite{FP94}.
    The preceding works \cite{KYN11} and \cite{YKN12} present formulae for computation of the diagonal entries of     $( B^{\top} B        )^{- M}$
                                                                                                                  and $( B        B^{\top} )^{- M}$ $( M = 1, 2, \dots )$.
    As the sum of these entries, the trace $J_{M} ( B )$ is obtained.
    These formulae are reviewed in Sections \ref{Rev_WS} and \ref{Rev_DiagType}.
    On the other hand, the derivation of the formula in the preceding work \cite{YKY14}
                       does not aim to obtain the diagonal entries of    $( B^{\top} B        )^{- M}$
                                                                      or $( B        B^{\top} )^{- M}$ $( M = 1, 2, \dots )$.
    This formula is reviewed in Section \ref{Rev_ResolType}.
    Note that we use the convention $\sum_{i = j}^{k} = 0$ if $j > k$ throughout this paper.

    \subsection{A review of the formula in the preceding work \cite{KYN11}}  \label{Rev_WS}

    In this subsection, we briefly review the formula for the trace $J_{M} ( B )$ for a fixed positive integer $M$ in the preceding work by Kimura {\it et al}. \cite{KYN11}.
    This formula gives the diagonal entries of the matrices $( B^{\top} B )^{- M}$ and $( B B^{\top} )^{- M}$ in a form of recurrence relation.
    As the sum of these diagonal entries, the trace $J_{M} ( B )$ is obtained.
    For a fixed positive integer $M$ and $m = 0, 1, \dots , M$,
    let the $( i, i )$-entry of the matrices $( B^{\top} B )^{- m}$ and $( B B^{\top} )^{- m}$ be denoted by $v_{i}^{( m )}$ and $w_{i}^{( m )}$, respectively.
    These diagonal entries are obtained by the following theorem.
    \begin{thm}  \label{Th_Diag_WS}
        Let $M$ be a fixed positive integer.
        All the diagonal entries $v_{i}^{(M)}$ and $w_{i}^{(M)}$ of inverse matrices $( B^{\top} B )^{- M}$ and $( B B^{\top} )^{- M}$, respectively, 
        are obtained through a finite number of arithmetics by using the following simple recurrence relations.
        The recurrence relations are
        \begin{align}
            &v_{i}^{( 0 )} =       1,                                                                                                 & &i = 1 , \dots , N,                        \\
            &w_{i}^{( 0 )} =       1,                                                                                                 & &i = 1 , \dots , N,                        \\
            &v_{N}^{( p )} = \frac{1}{q_{N}}                                                         w_{N    }^{( p - 1 )},           & &                       \label{recurr_N}   \\
            &v_{i}^{( p )} = \frac{1}{q_{i}} \left( e_{i    } v_{i + 1}^{( p )}  + z_{i}^{(p - 1)} - w_{i    }^{( p - 1 )} \right) ,  & &i = 1 , \dots , N - 1,                    \\
            &w_{1}^{( p )} = \frac{1}{q_{1}}                                                         v_{1    }^{( p - 1 )},           & &                                          \\
            &w_{i}^{( p )} = \frac{1}{q_{i}} \left( e_{i - 1} w_{i - 1}^{( p )}  + z_{i}^{(p - 1)} - v_{i    }^{( p - 1 )} \right) ,  & &i = 2 , \dots , N,     \label{recurr_w}   \\
            &z_{1}^{( q )} =       2                                                                 v_{1    }^{( q     )},           & &                       \label{recurr_z1}  \\  
            &z_{i}^{( q )} =                                  z_{i - 1}^{( q )}
                             +     2         \left(           v_{i    }^{( q )}                    - w_{i - 1}^{( q     )} \right) ,  & &i = 2 , \dots , N      \label{recurr_z} 
        \end{align}
        for integers $p$ and $q$ such that $1 \leq p \leq M$ and $0 \leq q \leq M - 1$.
        Instead of equations (\ref{recurr_z1}) and (\ref{recurr_z}), the following relations can be used:
        \begin{gather*}
            z_{N}^{( q )} = 2                            w_{N}^{( q )},                                                             \\  
            z_{i}^{( q )} = z_{i + 1}^{( q )} + 2 \left( w_{i}^{( q )} - v_{i + 1}^{( q )} \right) ,  \qquad i = 1 , \dots , N - 1
        \end{gather*}
        for integers $q$ such that $0 \leq q \leq M-1$.
    \end{thm}

    Substitution of $p = 1$ and $q = 0$ into the recurrence relations in Theorem \ref{Th_Diag_WS} gives the following remark.
    \begin{rmk}  \label{recurrence_M1}
        The recurrence relations from (\ref{recurr_N}) through (\ref{recurr_w}) in Theorem \ref{Th_Diag_WS} for $p = 1$ are simplified to the recurrence relations
        \begin{align}
            &v_{N}^{( 1 )} = \frac{1}{q_{N}},                                                   & &                        \label{vN1}  \\
            &v_{i}^{( 1 )} = \frac{1}{q_{i}} \left( e_{i    } v_{i + 1}^{( 1 )}  + 1 \right) ,  & &i = 1 , \dots , N - 1,  \label{vi1}  \\
            &w_{1}^{( 1 )} = \frac{1}{q_{1}},                                                   & &                        \label{w11}  \\
            &w_{i}^{( 1 )} = \frac{1}{q_{i}} \left( e_{i - 1} w_{i - 1}^{( 1 )}  + 1 \right) ,  & &i = 2 , \dots , N.      \label{wi1}
        \end{align}
        Theorem \ref{Th_Diag_WS} for $M = 1$ is reduced to these recurrence relations.
    \end{rmk}

    It is noted that there exist some preceding works on some limited cases in numerical analysis.
    \begin{rmk}  \label{PrecWorksNA}
        A formula related to the equations from (\ref{vN1}) to (\ref{wi1}) for computing the diagonal entries of the inverse $( B B^{\top} )^{-1}$ is known.
        See \cite{FP94, vonMatt97, Rutishauser90}, for example.
        Another formula for computing the diagonal entries of the inverse $( B B^{\top} )^{- 2}$ is presented by von Matt \cite{vonMatt97}.
    \end{rmk}

    \subsection{A review of the formula in the preceding work \cite{YKN12}}  \label{Rev_DiagType}

    In this subsection, we briefly review the preceding work by Yamashita, Kimura and Nakamura \cite{YKN12}.
    They derived a formula for computing the diagonal entries of the inverse powers     $( B^{\top} B        )^{- M}$
                                                                                    and $( B        B^{\top} )^{- M}$ $( M = 1, 2, \dots )$ in a form of recurrence relations
    starting from  the formula in \cite{KYN11} which is reviewed in the previous subsection.
    In contrast to the formula in \cite{KYN11} which includes subtraction in it in the cases of $M \geq 2$,
    the derived formula consists of only addition, multiplication and division among positive quantities.
    Namely, it is ``subtraction-free''.
    Any possibility of cancellation error in numerical computation is clearly excluded from this property.
    As the sum of these diagonal entries, the traces $J_{M} ( B )$ $( M = 1, 2, \dots )$ are obtained.

    For $p = 1, \dots , M$, let the $( i, i )$-entry of $( B^{\top} B )^{- p}$ and $( B B^{\top} )^{- p}$ be denoted by $v_{i}^{( p )}$ and $w_{i}^{( p )}$, respectively.
    Quantities $\check{B}_{i}$ $( i = 1, \dots , N )$, $F_{i}$ $( i = 1, \dots , N - 1 )$ and $\tilde{F}_{i}$ $( i = 2, \dots , N )$ are introduced as follows.
    \begin{dfn}  \label{defRatios}
        \begin{align}
            &\check{B}_{i} = \frac{1        }{q_{i}},                           & &i = 1, \dots , N,      \label{def_cB}  \\
            &       F _{i} = \frac{e_{i    }}{q_{i}}= e_{i    } \check{B}_{i},  & &i = 1, \dots , N - 1,  \label{def_F}   \\
            &\tilde{F}_{i} = \frac{e_{i - 1}}{q_{i}}= e_{i - 1} \check{B}_{i},  & &i = 2, \dots , N.      \label{def_tF}
        \end{align}
    \end{dfn}

    Next, quantities $g_{i}^{( r )}$ and $\tilde{g}_{i}^{( r )}$ for $i = 1, \dots , N $ and $r = 1, 2, \dots$ are introduced as follows.
    \begin{dfn}  \label{def_g}
        The quantities $g_{i}^{( r )}$ for $i = 1, \dots , N $ and $r = 1, 2, \dots$ are defined as follows.
        \begin{itemize}
            \item For $i =            N    $ and arbitrary positive integer $r$,               $g_{N}^{( r )}$  is given as $g_{N}^{( r )} = 0$.
            \item For $i = 1, \dots , N - 1$ and                            $r = 1$,           $g_{i}^{( 1 )}$  is given as $g_{i}^{( 1 )} = F_{i} v_{i + 1}^{( 1 )}$.
            \item For $i = 1, \dots , N - 1$ and                            $r = 2, 3, \dots$, $g_{i}^{( r )}$  is given as
                  \begin{equation}
                      g_{i}^{( r )} = F_{i} g_{i + 1}^{( r )} + \check{B}_{i + 1} g_{i}^{( r - 1 )} + \sum_{k = 1}^{r - 1} g_{i + 1}^{( k )} g_{i}^{( r - k )}.  \label{DefHigh_g}
                  \end{equation}
        \end{itemize}
    \end{dfn}
    \begin{dfn}  \label{def_tg}
        The quantities $\tilde{g}_{i}^{( r )}$ for $i = 1, \dots , N $ and $r = 1, 2, \dots$ are defined as follows.
        \begin{itemize}
            \item For $i =            1$ and arbitrary positive integer $r$,               $\tilde{g}_{1}^{( r )}$  is given as $\tilde{g}_{1}^{( r )} = 0$.
            \item For $i = 2, \dots , N$ and                            $r = 1$,           $\tilde{g}_{i}^{( 1 )}$  is given as $\tilde{g}_{i}^{( 1 )} = \tilde{F}_{i} w_{i - 1}^{( 1 )}$.
            \item For $i = 2, \dots , N$ and                            $r = 2, 3, \dots$, $\tilde{g}_{i}^{( r )}$  is given as
                  \begin{equation}  \label{DefHigh_tilde_g}
                      \tilde{g}_{i}^{( r )} =   \tilde{F}_{i    } \tilde{g}_{i - 1}^{( r     )}
                                              + \check{B}_{i - 1} \tilde{g}_{i    }^{( r - 1 )} + \sum_{k = 1}^{r - 1} \tilde{g}_{i - 1}^{( k )} \tilde{g}_{i}^{( r - k )}.
                  \end{equation}
        \end{itemize}
    \end{dfn}
    \begin{rmk}  \label{Rmk_M1_SF}
        The recurrence relations in Remark \ref{recurrence_M1} can be rewritten with the quantities defined by Definitions from \ref{defRatios} to \ref{def_tg} as follows.
        \begin{align}
            &v_{N}^{( 1 )} =                                         \check{B}_{N},  & &                       \label{vedge1}    \\
            &v_{i}^{( 1 )} =        F _{i}        v _{i + 1}^{(1)} + \check{B}_{i}
                           =                      g _{i    }^{(1)} + \check{B}_{i},  & &i = 1, \dots , N - 1,  \label{vrecurr1}  \\
            &w_{1}^{( 1 )} =                                         \check{B}_{1},  & &                       \label{wedge1}    \\
            &w_{i}^{( 1 )} = \tilde{F}_{i}        w _{i - 1}^{(1)} + \check{B}_{i}
                           =               \tilde{g}_{i    }^{(1)} + \check{B}_{i},  & &i = 2, \dots , N.      \label{wrecurr1}
        \end{align}
    \end{rmk}
    A formula for the traces $J_{M} ( B )$ $( M \geq 2 )$ is given as follows.
    \begin{thm}
        For $M \geq 2$, the diagonal entries $v_{i}^{( s )}$ and $w_{i}^{( s )}$ of $( B^{\top} B )^{- s}$ and $( B B^{\top} )^{- s}$, respectively,
        for $i = 1 , \dots , N$ and $s = 2 , \dots , M$ are computed by the recurrence relations
        \begin{align*}
            &v_{N}^{( s )} =                          \check{B}_{N}         w_{N    }^{(s - 1)},  & &                          \\
            &w_{1}^{( s )} =                          \check{B}_{1}         v_{1    }^{(s - 1)},  & &                          \\
            &v_{i}^{( s )} =                                 F _{i}         v_{i + 1}^{(s    )}
                             +                        \check{B}_{i}         w_{i    }^{(s - 1)}
                             + 2 \sum_{k = 1}^{s - 1}        g _{i}^{( k )} w_{i    }^{(s - k)},  & &  i = 1 , \dots , N - 1,  \\
            &w_{i}^{( s )} =                          \tilde{F}_{i}         w_{i - 1}^{(s    )}
                             +                        \check{B}_{i}         v_{i    }^{(s - 1)}
                             + 2 \sum_{k = 1}^{s - 1} \tilde{g}_{i}^{( k )} v_{i    }^{(s - k)},  & &  i = 2 , \dots , N,
        \end{align*}
        with the recurrence relations from (\ref{vedge1}) to (\ref{wrecurr1}).
    \end{thm}

    \subsection{A review of the formula in the preceding work \cite{YKY14}}  \label{Rev_ResolType}

    In this subsection, we briefly review the preceding work by Yamashita, Kimura and Yamamoto \cite{YKY14}.
    Since the derivation of the formula does not aim to obtain the diagonal entries of $( B^{\top} B )^{- p}$ or $( B B^{\top} )^{- p}$ $( p = 1, 2, \dots )$,
    we review the idea for the derivation briefly first.
    Let us consider the matrix $B^{\top} B$.
    It is transformed into
    \begin{equation*}
        A = \left( \begin{array}{cccc}
                       q_{1}  &  q_{1}   e_{1}  &          &                         \\
                       1      &  q_{2} + e_{1}  &  \ddots  &                         \\
                              &  \ddots         &  \ddots  &  q_{N - 1}   e_{N - 1}  \\
                              &                 &  1       &  q_{N    } + e_{N - 1}
                   \end{array}                                                           \right)
    \end{equation*}
    by similarity transformation without changing its eigenvalues.
    Let $I$ be the $N \times N$ unit matrix.
    In derivation of the formula from the matrix $B^{\top} B$,
    one of the key points is to express the determinant $D ( \lambda ) = \det ( A - \lambda I )$ in two ways.
    For the eigenvalues $\lambda _{1}, \dots , \lambda _{N}$ of $B^{\top} B$,
    the determinant is given as $D ( \lambda ) = \prod_{i = 1}^{N} ( \lambda _{i} - \lambda )$.
    The matrix $A - \lambda I$ can be decomposed as
    \begin{equation*}
        A - \lambda I = \left( \begin{array}{cccc}
                                    \hat{q}_{1}^{( 0 )}  &                       &          &                       \\
                                    1                    &  \hat{q}_{2}^{( 0 )}  &          &                       \\
                                                         &  \ddots               &  \ddots  &                       \\
                                                         &                       &  1       &  \hat{q}_{N}^{( 0 )}
                                \end{array}                                                                            \right)
                         \left( \begin{array}{cccc}
                                    1  &  \hat{e}_{1}^{( 0 )}  &          &                           \\
                                       &  1                    &  \ddots  &                           \\
                                       &                       &  \ddots  &  \hat{e}_{N - 1}^{( 0 )}  \\
                                       &                       &          &  1
                                \end{array}                                                                            \right) ,
    \end{equation*}
    where $\hat{q}_{i}^{( 0 )}$ for $i = 1, \dots , N$ and $\hat{e}_{i}^{( 0 )}$ for $i = 1, \dots , N - 1$ are functions of $\lambda$.
    The superscript $( 0 )$ means that these functions have not been differentiated yet.
    The other expression of the determinant is $D ( \lambda ) = \prod_{i = 1}^{N} \hat{q}_{i}^{( 0 )}$.
    Since the matrices $B^{\top} B$ and $A$ have the same eigenvalues, it holds
    \begin{equation}  \label{EqDetAlI}
        \prod_{i = 1}^{N} ( \lambda _{i} - \lambda ) = \prod_{i = 1}^{N} \hat{q}_{i}^{( 0 )}.
    \end{equation}
    Let $\hat{q}_{i}^{( p )}$ for $i = 1, \dots , N$ and $p = 1, 2, \dots$ be $d^{p} \hat{q}_{i}^{( 0 )} / d \lambda ^{p}$.
    By differentiating (\ref{EqDetAlI}) with respect to $\lambda$, we have
    \begin{equation*}
          \sum_{i = 1}^{N} \left( - \frac{D ( \lambda )      }{\lambda _{i}         - \lambda}                     \right)
        = \sum_{i = 1}^{N} \left(   \frac{\hat{q}_{i}^{( 1 )}}{\hat{q} _{i}^{( 0 )}          } \cdot D ( \lambda ) \right) .
    \end{equation*}
    Then, it holds
    \begin{equation}  \label{DiffOnce}
        \sum_{i = 1}^{N} \frac{1}{\lambda _{i} - \lambda} = \sum_{i = 1}^{N} \left( - \frac{\hat{q}_{i}^{( 1 )}}{\hat{q}_{i}^{( 0 )}} \right) .
    \end{equation}
    Let $\hat{H}_{i}^{( p )}$ for $i = 1, \dots , N$ and $p = 1, 2, \dots$ be functions of $\lambda$ defined as 
    \begin{equation}
        \hat{H}_{i}^{( p )} = \left \{ \begin{array}{ll}
                                  \displaystyle - \frac{  \hat{q}_{i}^{(     1 )}}{  \hat{q}_{i}^{( 0 )}},  &  p = 1,             \\
                                  \displaystyle   \frac{d \hat{H}_{i}^{( p - 1 )}}{d \lambda            },  &  p = 2, 3, \dots .
                              \end{array} \right .
    \end{equation}
    Then, by differentiating (\ref{DiffOnce}) repeatedly, we have
    \begin{equation}  \label{lambda_H}
          ( p - 1 )! \sum_{i = 1}^{N} \frac{1}{( \lambda _{i} - \lambda )^{p}}
        =            \sum_{i = 1}^{N} \hat{H}_{i}^{( p )},                      \qquad  p = 1, 2, \dots .
    \end{equation}
    Thus, the traces of $( B^{\top} B )^{- p}$ $( p = 1, 2, \dots )$ are obtained by substituting $\lambda = 0$ into (\ref{lambda_H}).
    The result is
    \begin{equation}  \label{TrH_all}
        \textrm{Tr} ( ( B^{\top} B )^{- p} ) = \frac{1}{( p - 1 )!} \sum_{i = 1}^{N} H_{i}^{( p )},  \qquad  p = 1, 2, \dots ,
    \end{equation}
    where $H_{i}^{( p )}$ is $H_{i}^{( p )} = \left. \hat{H}_{i}^{( p )} \right| _{\lambda = 0}$.
    From this result, a method for computing $H_{i}^{( p )}$ is required.
    As such a method, a recurrence relation with a help of auxiliary quantities $h_{i}^{( p )}$ for $i = 1, \dots , N$ and $p = 1, 2, \dots$ is given.
    The auxiliary quantities $h_{i}^{( p )}$ is defined as $h_{i}^{( p )} = \left. \hat{h}_{i}^{( p )} \right| _{\lambda = 0}$,
    where $\hat{h}_{i}^{( p )}$ are functions of $\lambda$ defined as $\hat{h}_{i}^{( p )} = - \hat{q}_{i}^{( p )} / \hat{q}_{i}^{( 0 )}$.
    The recurrence relation is presented as the following theorem.
    \begin{thm} \label{Th_G}
        Let $B$ be an upper bidiagonal matrix defined in (\ref{UpperBidiagonal}).

        Let us introduce constants $\tilde{F}_{i}$ for $i = 2, \dots , N$ defined as
        \begin{equation*}
            \tilde{F}_{i} = \frac{e_{i - 1}}{q_{i}},  \qquad  i = 2, \dots , N.
        \end{equation*}
        For $i = 1, \dots , N$ and $p = 1, 2, \dots$, let $h_{i}^{( p )}$ be constants which satisfy the following recurrence relation
        \begin{align}
            &h_{1}^{( 1 )} =                                   \frac{1}{q_{1}},                                        & &                   & &                   \label{h11}  \\
            &h_{i}^{( 1 )} = \tilde{F}_{i} h_{i - 1}^{( 1 )} + \frac{1}{q_{i}},                                        & &i = 2, \dots , N,  & &                   \label{hi1}  \\
            &h_{1}^{( p )} = 0,                                                                                        & &                   & &p = 2, 3, \dots ,  \label{h1h}  \\
            &h_{i}^{( p )} = \tilde{F}_{i} \left(                                       h_{i - 1}^{( p     )}
                                                  + p                                   h_{i - 1}^{(     1 )}
                                                                                        h_{i - 1}^{( p - 1 )} \right)  & &                   & &                   \notag       \\
            & \qquad                              +   \sum_{k = 1}^{p - 2} {}_{p} C_{k} h_{i - 1}^{( k     )}
                                                                                        h_{i    }^{( p - k )},         & &i = 2, \dots , N,  & &p = 2, 3, \dots .  \label{hih}
        \end{align}
        For $i = 1, \dots , N$ and $p = 1, 2, \dots$, let $H_{i}^{( p )}$ be constants which satisfy the following recurrence relation
        \begin{align}
            &H_{i}^{( 1 )} = h_{i}^{( 1 )},                                                                          & &i = 1, \dots , N,  & &                   \notag         \\
            &H_{i}^{( p )} = h_{i}^{( p )} + \sum_{k = 1}^{p - 1} {}_{p - 1} C_{k} h_{i}^{( k )} H_{i}^{( p - k )},  & &i = 1, \dots , N,  & &p = 2, 3, \dots .  \label{Hhigh}
        \end{align}

        The traces $\textrm{Tr} ( ( B^{\top} B )^{- p} )$ for $p = 1, 2, \dots$ are computed by
        \begin{equation*}
            \textrm{Tr} ( ( B^{\top} B )^{- p} ) = \frac{1}{( p - 1 ) !} \sum_{i = 1}^{N} H_{i}^{( p )}.
        \end{equation*}
    \end{thm}
    \begin{rmk}  \label{Rmk_diff_w}
        The recurrence relation of $h_{i}^{( 1 )}$ $( i = 1, \dots , N )$ in Theorem \ref{Th_G} is equivalent
        to the recurrence relation for the diagonal entries of $( B B^{\top} )^{- 1}$ shown in Remark \ref{recurrence_M1} (Originally, Remark 4.6 in \cite{KYN11}.).
        Then, the constants $h_{i}^{( 1 )}$ and $H_{i}^{( 1 )}$ are the $( i, i )$-entry of $( B B^{\top} )^{- 1}$.
        See also Remark \ref{PrecWorksNA} (Originally, Remark 4.7 in \cite{KYN11}.).
    \end{rmk}

    Considering the matrix $B B^{\top}$ instead of $B^{\top} B$, we have another recurrence relation as follows.
    \begin{thm} \label{Th_tG}
        Let $B$ be an upper bidiagonal matrix defined in (\ref{UpperBidiagonal}).

        Let us introduce constants $F_{i}$ for $i = 1, \dots , N - 1$ defined as
        \begin{equation*}
            F_{i} = \frac{e_{i}}{q_{i}},  \qquad  i = 1, \dots , N - 1.
        \end{equation*}
        For $i = 1, \dots , N$ and $p = 1, 2, \dots$, let $\tilde{h}_{i}^{( p )}$ be constants which satisfy the following recurrence relation
        \begin{align*}
            &\tilde{h}_{N}^{( 1 )} =                                   \frac{1}{q_{N}},                                        & &                       & &                   \\
            &\tilde{h}_{i}^{( 1 )} = F_{i} \tilde{h}_{i + 1}^{( 1 )} + \frac{1}{q_{i}},                                        & &i = 1, \dots , N - 1,  & &                   \\
            &\tilde{h}_{N}^{( p )} = 0,                                                                                        & &                       & &p = 2, 3, \dots ,  \\
            &\tilde{h}_{i}^{( p )} = F_{i} \left(                                       \tilde{h}_{i + 1}^{( p     )}
                                                  + p                                   \tilde{h}_{i + 1}^{(     1 )}
                                                                                        \tilde{h}_{i + 1}^{( p - 1 )} \right)  & &                       & &                   \\
            & \qquad                              +   \sum_{k = 1}^{p - 2} {}_{p} C_{k} \tilde{h}_{i + 1}^{( k     )}
                                                                                        \tilde{h}_{i    }^{( p - k )},         & &i = 1, \dots , N - 1,  & &p = 2, 3, \dots .
        \end{align*}
        For $i = 1, \dots , N$ and $p = 1, 2, \dots$, let $\tilde{H}_{i}^{( p )}$ be constants which satisfy the following recurrence relation
        \begin{align*}
            &\tilde{H}_{i}^{( 1 )} =                                         \tilde{h}_{i}^{( 1 )},                            & &i = 1, \dots , N,  & &                   \\
            &\tilde{H}_{i}^{( p )} =                                         \tilde{h}_{i}^{( p )}
                                     + \sum_{k = 1}^{p - 1} {}_{p - 1} C_{k} \tilde{h}_{i}^{( k )} \tilde{H}_{i}^{( p - k )},  & &i = 1, \dots , N,  & &p = 2, 3, \dots .
        \end{align*}

        The traces $\textrm{Tr} ( ( B B^{\top} )^{- p} )$ for $p = 1, 2, \dots$ are computed by
        \begin{equation}  \label{TraceByH}
            \textrm{Tr} ( ( B B^{\top} )^{- p} ) = \frac{1}{( p - 1 ) !} \sum_{i = 1}^{N} \tilde{H}_{i}^{( p )}.
        \end{equation}
    \end{thm}
    \begin{rmk}
        The recurrence relation of $\tilde{h}_{i}^{( 1 )}$ $( i = 1, \dots , N )$ in Theorem \ref{Th_tG} is equivalent
        to the recurrence relation for the diagonal entries of $( B^{\top} B )^{- 1}$ shown in Remark \ref{recurrence_M1} (Originally, Remark 4.6 in \cite{KYN11}.).
        Then, the constants $\tilde{h}_{i}^{( 1 )}$ and $\tilde{H}_{i}^{( 1 )}$ are the $( i, i )$-entry of $( B^{\top} B )^{- 1}$.
    \end{rmk}

    The recurrence relations presented in Theorems \ref{Th_G} and \ref{Th_tG} are subtraction-free.

\section{A unified understanding of the two formulae in the preceding works and a new subtraction-free formula for the traces}  \label{NewFormula}

    In this section, we achieve a unified understanding of the two formulae in the preceding works \cite{YKN12} and \cite{YKY14} which are reviewed in the previous section.
    Firstly, we give an interpretation of      some quantities, which is introduced in the formula  in \cite{YKN12}, in terms of matrix theory.
    Next,    we show a  relationship   between some quantities                      in the formulae in \cite{YKN12} and \cite{YKY14}.
    Lastly,  a new subtraction-free formula for the traces $J_{M} ( B )$ $( M = 1, 2, \dots )$ is naturally derived from this relationship.
    
    Let us consider the quantities $\tilde{g}_{i}^{( M )}$ for $i = 1, \dots , N$ and $M = 1, 2, \dots$ given in Definition \ref{def_tg}.
    We give an interpretation of the quantities $\tilde{g}_{i}^{( M )}$ for $i = 1, \dots , N$ and $M = 2, 3, \dots$ in terms of matrix theory.
    The following theorem holds.
    \begin{thm}  \label{Th_g_Mat}
        Let $B$ be an upper bidiagonal matrix defined in (\ref{UpperBidiagonal}).
        Let $W = ( W_{i, j} )$ denote the matrix $( B B^{\top} )^{- 1}$.
        The quantities $\tilde{g}_{i}^{( M )}$ for $i = 1, \dots , N$ and $M = 2, 3, \dots$ which are given in Definition \ref{def_tg} are expressed as
        \begin{gather*}
            \tilde{g}_{i}^{( 2 )} = \sum_{j_{1    } = 1}^{i - 1} W_{i, j_{1}}                                                  W_{j_{1    }, i},                             \\
            \tilde{g}_{i}^{( 3 )} = \sum_{j_{1    } = 1}^{i - 1}
                                    \sum_{j_{2    } = 1}^{i - 1} W_{i, j_{1}} W_{j_{1}, j_{2}}                                 W_{j_{2    }, i},                             \\
            \tilde{g}_{i}^{( M )} = \sum_{j_{1    } = 1}^{i - 1}
                                    \cdots
                                    \sum_{j_{M - 1} = 1}^{i - 1} W_{i, j_{1}} W_{j_{1}, j_{2}} \cdots W_{j_{M - 2}, j_{M - 1}} W_{j_{M - 1}, i},  \qquad  M = 4, 5, \dots .
        \end{gather*}
    \end{thm}

    Proof of this theorem is given in Appendix.

    Next, we reveal a relationship between the quantities $\tilde{g}_{i}^{( M )}$ and $h_{i}^{( M )}$ for $i = 1, \dots , N$ and $M = 2, 3, \dots$.
    As a preliminary, let us derive expressions for the traces $J_{2} ( B )$ and $J_{3} ( B )$
                                    with $\tilde{g}_{i}^{( 2 )}$, $\tilde{g}_{i}^{( 3 )}$ and $W_{i, i}( = h_{i}^{( 1 )} )$ $( i = 1, \dots , N )$.
    The trace $J_{2} ( B )$ is given as
    \begin{equation}  \label{Tr2W}
        J_{2} ( B ) = \sum_{i = 1}^{N} \sum_{j = 1}^{N} W_{i, j} W_{j, i}.
    \end{equation}
    Let us classify the products $W_{i, j} W_{j, i}$ in (\ref{Tr2W}) into three sets according to the subscripts $i$ and $j$.
    Namely, i) $i > j$, ii) $i < j$ and iii) $i = j$.
    The summation of the products which belong to the first set is
    \begin{equation*}
        \sum_{i = 1}^{N} \sum_{j = 1}^{i - 1} W_{i, j} W_{j, i} = \sum_{i = 1}^{N} \tilde{g}_{i}^{( 2 )}
    \end{equation*}
    from Theorem \ref{Th_g_Mat}.
    That of the second set is equal to this summation.
    That of the third  set is
    \begin{equation*}
        \sum_{i = 1}^{N} \left( W_{i, i} \right) ^{2} = \sum_{i = 1}^{N} \left( h_{i}^{( 1 )} \right) ^{2}.
    \end{equation*}
    Then, the trace $J_{2} ( B )$ is written as
    \begin{equation}  \label{Tr2g}
        J_{2} ( B ) = \sum_{i = 1}^{N} \left( 2 \tilde{g}_{i}^{( 2 )} + \left( h_{i}^{( 1 )} \right) ^{2} \right) .
    \end{equation}
    The trace $J_{3} ( B )$ is given as
    \begin{equation}  \label{Tr3W}
        J_{3} ( B ) = \sum_{i = 1}^{N} \sum_{j = 1}^{N} \sum_{k = 1}^{N} W_{i, j} W_{j, k} W_{k, i}.
    \end{equation}
    The products $W_{i, j} W_{j, k} W_{k, i}$ in (\ref{Tr3W}) are classified into seven sets according to the subscripts $i$, $j$ and $k$.
    Namely, i) $i > j$ and $i > k$, ii) $j > i$ and $j > k$, iii) $k > i$ and $k > j$, iv) $i = j > k$, v) $i = k > j$, vi) $j = k > i$ and vii) $i = j = k$.
    The summation of the products which belong to the first set is
    \begin{equation*}
        \sum_{i = 1}^{N} \sum_{j = 1}^{i - 1} \sum_{k = 1}^{i - 1}  W_{i, j} W_{j, k} W_{k, i} = \sum_{i = 1}^{N} \tilde{g}_{i}^{( 3 )}
    \end{equation*}
    from Theorem \ref{Th_g_Mat}.
    Those of the second and the third sets are equal to this summation.
    For the fourth set, the summation is
    \begin{equation*}
        \sum_{i = 1}^{N} \sum_{k = 1}^{i - 1}  W_{i, i} W_{i, k} W_{k, i} = \sum_{i = 1}^{N} h_{i}^{( 1 )} \tilde{g}_{i}^{( 2 )}
    \end{equation*}
    from Theorem \ref{Th_g_Mat}.
    Those of the fifth and the sixth sets are equal to this summation.
    For the seventh set, the summation is
    \begin{equation*}
        \sum_{i = 1}^{N} \left( W_{i, i} \right) ^{3} = \sum_{i = 1}^{N} \left( h_{i}^{( 1 )} \right) ^{3}.
    \end{equation*}
    Thus, the trace $J_{3} ( B )$ is expressed as
    \begin{equation}  \label{Tr3g}
        J_{3} ( B ) = \sum_{i = 1}^{N} \left( 3 \tilde{g}_{i}^{( 3 )} + 3 \tilde{g}_{i}^{( 2 )} h_{i}^{( 1 )} + \left( h_{i}^{( 1 )} \right) ^{3} \right) .
    \end{equation}
    On the other hand, it is possible to express these traces with $h_{i}^{( 1 )}$, $h_{i}^{( 2 )}$ and $h_{i}^{( 3 )}$ $( i = 1, \dots , N )$.
    From Theorem \ref{Th_G}, we readily have
    \begin{gather}
        J_{2} ( B ) =             \sum_{i = 1}^{N} \left( h_{i}^{( 2 )} +                 \left( h_{i}^{( 1 )} \right) ^{2} \right) ,  \label{Tr2h}  \\
        J_{3} ( B ) = \frac{1}{2} \sum_{i = 1}^{N} \left( h_{i}^{( 3 )} + 3 h_{i}^{( 2 )}        h_{i}^{( 1 )}
                                                                        + 2               \left( h_{i}^{( 1 )} \right) ^{3} \right) .  \label{Tr3h}
    \end{gather}
    Comparison of (\ref{Tr2g}) with (\ref{Tr2h}) and of (\ref{Tr3g}) with (\ref{Tr3h}) suggests us to investigate whether it holds that
    \begin{equation*}
        h_{i}^{( M )} = M ! \cdot \tilde{g}_{i}^{( M )},  \qquad  i = 1, \dots , N,  \qquad  M = 2, 3, \dots ,
    \end{equation*}
    or not.
    The answer of the investigation is ``true''.
    The following theorem holds.
    \begin{thm}  \label{Th_trans_gh}
        For $i = 1, \dots , N$ and $M = 2, 3, \dots$, the transform between a quantities $\tilde{g}_{i}^{( M )}$ and $h_{i}^{( M )}$ is given as
        \begin{equation}  \label{trans_gh}
             h_{i}^{( M )} = M ! \cdot \tilde{g}_{i}^{( M )}.
        \end{equation}
    \end{thm}

    Proof of this theorem is given in Appendix.

    Now, we derive the new subtraction-free formula.
    Let us introduce the following transform
    \begin{equation}  \label{Transf_HtG}
        H_{i}^{( M )} = ( M - 1 ) ! \cdot \tilde{G}_{i}^{( M )},  \qquad  i = 1, \dots , N,  \qquad  M = 1, 2, \dots .
    \end{equation}
    Then, we readily have a new formula for the traces $J_{M} ( B )$ $( M = 1, 2, \dots )$ as follows by substituting (\ref{trans_gh}) and (\ref{Transf_HtG})
    into the recurrence relation in Theorem \ref{Th_G} in Section \ref{Rev_ResolType}.
    \begin{thm}  \label{Th_NewFormula}
        Let $B$ be an $N \times N$ upper bidiagonal matrix defined as (\ref{UpperBidiagonal}).
        Let $\check{B}_{i}$ for $i = 1, \dots , N$ be quantities defined by (\ref{def_cB}).
        Let $\tilde{F}_{i}$ for $i = 2, \dots , N$ be quantities defined by (\ref{def_tF}).

        Let $\tilde{g}_{i}^{( 1 )}$ and $\tilde{G}_{i}^{( 1 )}$ for $i = 1, \dots , N$ be quantities which satisfy the recurrence relation
        \begin{align}
            &\tilde{g}_{1}^{( 1 )} = 0,                                                        & &                   \label{tg11}  \\
            &\tilde{G}_{1}^{( 1 )} =                                           \check{B}_{1},  & &                   \label{tG11}  \\
            &\tilde{g}_{i}^{( 1 )} = \tilde{F}_{i} \tilde{G}_{i - 1}^{( 1 )},                  & &i = 2, \dots , N,  \label{tgi1}  \\
            &\tilde{G}_{i}^{( 1 )} =               \tilde{g}_{i    }^{( 1 )} + \check{B}_{i},  & &i = 2, \dots , N.  \label{tGi1}
        \end{align}
        Let $\tilde{g}_{i}^{( M )}$ for $i = 1, \dots , N$ and $M = 2, 3, \dots$ be quantities which satisfy the following recurrence relation
        \begin{gather}
            \tilde{g}_{1}^{( M )} = 0,                                                                                                          \label{tg1h} \\
            \tilde{g}_{i}^{( M )} =                        \tilde{F}_{i    }         \tilde{g}_{i - 1}^{( M     )}
                                    +                      \tilde{G}_{i - 1}^{( 1 )} \tilde{g}_{i    }^{( M - 1 )}
                                    + \sum_{k = 2}^{M - 1} \tilde{g}_{i - 1}^{( k )} \tilde{g}_{i    }^{( M - k )},  \qquad  i = 2, \dots , N.  \label{tgih}
        \end{gather}
        Let $\tilde{G}_{i}^{( M )}$ for $i = 1, \dots , N$ and $M = 2, 3, \dots$ be quantities which satisfy the following recurrence relation
        \begin{equation}  \label{tGhigh}
            \tilde{G}_{i}^{( M )} =   M                      \tilde{g}_{i}^{( M )}
                                    +                        \tilde{G}_{i}^{( 1 )} \tilde{G}_{i}^{( M - 1 )}
                                    +   \sum_{k = 2}^{M - 1} \tilde{g}_{i}^{( k )} \tilde{G}_{i}^{( M - k )}.
        \end{equation}

        Then, the trace $J_{M} ( B )$ for an arbitrary positive integer $M$ is obtained from
        \begin{equation}  \label{TraceByTG}
            J_{M} ( B ) = \sum_{i = 1}^{N} \tilde{G}_{i}^{( M )}.
        \end{equation}
    \end{thm}

    Proof of this theorem is given in Appendix.

    As a comparison of the formulae in Theorems \ref{Th_G} and \ref{Th_NewFormula}, we have the following remark.
    \begin{rmk}
        In contrast to the formula in Theorem \ref{Th_G}, the formula in Theorem \ref{Th_NewFormula} includes neither combination nor factorial in it.
        The quantities $\tilde{G}_{i}^{( M )}$ $( i = 1, \dots , N )$ in Theorem \ref{Th_NewFormula} directly gives the trace $J_{M} ( B )$ as their summation 
        in contrast to the quantities $H_{i}^{( M )}$ $( i = 1, \dots , N )$ in Theorem \ref{Th_G}.
    \end{rmk}
    
    The relationship (\ref{Transf_HtG}) for $M = 1$ and Remark \ref{Rmk_diff_w} give the following remark.
    \begin{rmk}  \label{Rmk_diag_NewFormula}
        The quantities $\tilde{G}_{i}^{( 1 )}$ $( i = 1, \dots , N )$ are the $( i, i )$-entry of $( B B^{\top} )^{- 1}$.
    \end{rmk}

    Let us consider the quantities $g_{i}^{( M )}$ for $i = 1, \dots , N$ and $M = 1, 2, \dots$ given in Definition \ref{def_g} instead of $\tilde{g}_{i}^{( M )}$.
    Similarly to above, we have the following theorems and remarks.
    \begin{thm}  \label{Th_tg_Mat}
        Let $B$ be an upper bidiagonal matrix defined in (\ref{UpperBidiagonal}).
        Let $V = ( V_{i, j} )$ denote the matrix $( B^{\top} B )^{- 1}$.
        The quantities $g_{i}^{( M )}$ for $i = 1, \dots , N$ and $M = 2, 3, \dots$ which are given in Definition \ref{def_g} are expressed as
        \begin{gather*}
            g_{i}^{( 2 )} = \sum_{j_{1    } = i + 1}^{N} V_{i, j_{1}}                                                  V_{j_{1    }, i},                             \\
            g_{i}^{( 3 )} = \sum_{j_{1    } = i + 1}^{N}
                            \sum_{j_{2    } = i + 1}^{N} V_{i, j_{1}} V_{j_{1}, j_{2}}                                 V_{j_{2    }, i},                             \\
            g_{i}^{( M )} = \sum_{j_{1    } = i + 1}^{N}
                            \cdots
                            \sum_{j_{M - 1} = i + 1}^{N} V_{i, j_{1}} V_{j_{1}, j_{2}} \cdots V_{j_{M - 2}, j_{M - 1}} V_{j_{M - 1}, i},  \qquad  M = 4, 5, \dots .
        \end{gather*}
    \end{thm}

    \begin{thm}  \label{Th_trans_gh_pair}
        For $i = 1, \dots , N$ and $M = 2, 3, \dots$, a transform between the quantities $g_{i}^{( M )}$ and $\tilde{h}_{i}^{( M )}$ is given as
        \begin{equation*}
             \tilde{h}_{i}^{( M )} = M ! \cdot g_{i}^{( M )}.
        \end{equation*}
    \end{thm}

    \begin{thm}  \label{Th_NewFormula_pair}
        Let $B$ be an $N \times N$ upper bidiagonal matrix defined as (\ref{UpperBidiagonal}).
        Let $\check{B}_{i}$ for $i = 1, \dots , N    $ be quantities defined by (\ref{def_cB}).
        Let $       F _{i}$ for $i = 1, \dots , N - 1$ be quantities defined by (\ref{def_F}).

        Let $g_{i}^{( 1 )}$ and $G_{i}^{( 1 )}$ for $i = 1, \dots , N$ be quantities which satisfy the recurrence relation
        \begin{align*}
            &g_{N}^{( 1 )} = 0,                                        & &                       \\
            &G_{N}^{( 1 )} =                           \check{B}_{N},  & &                       \\
            &g_{i}^{( 1 )} = F_{i} G_{i + 1}^{( 1 )},                  & &i = 1, \dots , N - 1,  \\
            &G_{i}^{( 1 )} =       g_{i    }^{( 1 )} + \check{B}_{i},  & &i = 1, \dots , N - 1.
        \end{align*}
        Let $g_{i}^{( M )}$ for $i = 1, \dots , N$ and $M = 2, 3, \dots$ be quantities which satisfy the following recurrence relation
        \begin{gather*}
            g_{N}^{( M )} = 0,                                                                                          \\
            g_{i}^{( M )} =                        F_{i    }         g_{i + 1}^{( M     )}
                            +                      G_{i + 1}^{( 1 )} g_{i    }^{( M - 1 )}
                            + \sum_{k = 2}^{M - 1} g_{i + 1}^{( k )} g_{i    }^{( M - k )},  \qquad  i = 2, \dots , N.
        \end{gather*}
        Let $G_{i}^{( M )}$ for $i = 1, \dots , N$ and $M = 2, 3, \dots$ be quantities which satisfy the following recurrence relation
        \begin{equation*}
            G_{i}^{( M )} =   M                      g_{i}^{( M )}
                            +                        G_{i}^{( 1 )} G_{i}^{( M - 1 )}
                            +   \sum_{k = 2}^{M - 1} g_{i}^{( k )} G_{i}^{( M - k )}.
        \end{equation*}

        Then, the trace $J_{M} ( B )$ for an arbitrary positive integer $M$ is obtained from
        \begin{equation*}
            J_{M} ( B ) = \sum_{i = 1}^{N} G_{i}^{( M )}.
        \end{equation*}
    \end{thm}
    \begin{rmk}
        In contrast to the formula in Theorem \ref{Th_tG}, the formula in Theorem \ref{Th_NewFormula_pair} includes neither combination nor factorial in it.
        The quantities $G_{i}^{( M )}$ $( i = 1, \dots , N )$ in Theorem \ref{Th_NewFormula_pair} directly gives the trace $J_{M} ( B )$ as their summation 
        in contrast to the quantities $\tilde{H}_{i}^{( M )}$ $( i = 1, \dots , N )$ in Theorem \ref{Th_tG}.
    \end{rmk}
    \begin{rmk}
        The quantities $G_{i}^{( 1 )}$ $( i = 1, \dots , N )$ are the $( i, i )$-entry of $( B^{\top} B )^{- 1}$.
    \end{rmk}

\section{Concluding Remarks}  \label{ConcRmks}

    In this paper, we achieve a unified understanding of the formulae
    for the traces of the inverse powers of a positive definite symmetric tridiagonal matrix in the two preceding works.
    We show that the quantities $g_{i}^{( r )}$ and $\tilde{g}_{i}^{( r )}$ which are introduced in \cite{YKN12} and are shown in Section \ref{Rev_DiagType}
                 can be interpreted in terms of matrix theory.
    The relationships which hold between some quantities in the two formulae in \cite{YKN12} and \cite{YKY14} are revealed.
    From these relationships, a new subtraction-free formula for the traces is obtained.

\section*{Appendix}

    In this appendix, we give proof of Theorems \ref{Th_g_Mat}, \ref{Th_trans_gh} and \ref{Th_NewFormula}.
    For proof of Theorems \ref{Th_g_Mat} and \ref{Th_trans_gh}, we give the following remark.
    \begin{rmk}  \label{Rmk_recur_g}
        For $i = 2, \dots , N$, the quantities $\tilde{g}_{i}^{( 1 )}$ can be written as
        \begin{equation}  \label{rew_g1}
            \tilde{g}_{i}^{( 1 )} = \tilde{F}_{i} W_{i - 1, i - 1}
        \end{equation}
        from its definition.
        For $i = 2, \dots , N$, we can readily verify
                                that it holds $\tilde{g}_{i - 1}^{( 1 )} + \check{B}_{i - 1} = W_{i - 1, i - 1}$ from Definition \ref{def_tg} and Remark \ref{Rmk_M1_SF}.
        Then, the recurrence relation (\ref{DefHigh_tilde_g}) in Definition \ref{def_tg} for $i = 2, \dots , N$ and $M = 2, 3, \dots$ can be rewritten as
        \begin{equation}  \label{rew_gh}
            \tilde{g}_{i}^{( M )} =                        \tilde{F}_{i    }                          \tilde{g}_{i - 1}^{( M     )}
                                    +                                                W_{i - 1, i - 1} \tilde{g}_{i    }^{( M - 1 )}
                                    + \sum_{k = 2}^{M - 1} \tilde{g}_{i - 1}^{( k )}                  \tilde{g}_{i    }^{( M - k )}.
        \end{equation}
    \end{rmk}

    Before showing proof of Theorem \ref{Th_g_Mat}, we give a lemma for this proof.
    \begin{lem}  \label{LemmaEntW}
        If $i > j$ and $i > k$, then it holds
        \begin{gather}
            W_{i, j    } W_{k, i} = \tilde{F}_{i} W_{i - 1, j} W_{k, i - 1},  \label{WW_dd}  \\
            W_{j, i - 1} W_{k, i} =               W_{j    , i} W_{k, i - 1}.  \label{WW_ud}
        \end{gather}
    \end{lem}
    {\it proof}

    Let $S = ( S_{i, j} )$ denote the inverse of $B$.
    Then, the matrix $W$ is expressed as $W = S^{\top} S$ from the definition $W = ( B B^{\top} )^{- 1}$.
    This matrix $S$ is an upper triangle matrix.
    The relationship
    \begin{equation*}
        S_{i, j} = - \sqrt{\tilde{F}_{j}} \cdot S_{i, j - 1},  \qquad  1 \leq i < j \leq N
    \end{equation*}
    holds among entries of $S$ (See \cite{KYN11}.).
    This relationship gives relationships among entries of $W$.
    When $i > j$, since it holds
    \begin{align*}
        W_{i, j} &=                        \sum_{k = 1}^{N}      (                              S^{\top} )_{i    , k    }         S_{k, j}
                  =                        \sum_{k = 1}^{N}                                     S         _{k    , i    }         S_{k, j}
                  =                        \sum_{k = 1}^{j}                                     S         _{k    , i    }         S_{k, j}  \\
                 &=                        \sum_{k = 1}^{j} \left( - \sqrt{\tilde{F}_{i}} \cdot S         _{k    , i - 1} \right) S_{k, j}
                  = - \sqrt{\tilde{F}_{i}} \sum_{k = 1}^{N}      (                              S^{\top} )_{i - 1, k    }         S_{k, j}  \\
                 &= - \sqrt{\tilde{F}_{i}}                                                \cdot W         _{i - 1, j    },
    \end{align*}
    we have
    \begin{equation*}
        W_{i, j} = - \sqrt{\tilde{F}_{i}} \cdot W_{i - 1, j}.
    \end{equation*}
    Since $W$ is a symmetric matrix, it holds
    \begin{gather*}
        W_{j, i} = - \sqrt{\tilde{F}_{i}} \cdot W_{j, i - 1}
    \end{gather*}
    when $i > j$.
    Then, when $i > j$ and $i > k$, it holds (\ref{WW_dd}) and (\ref{WW_ud}).  ~$\square$

\subsection*{Proof of Theorem \ref{Th_g_Mat}}
    
    Let us introduce quantities $u_{i}^{( M )}$ for $i = 1, \dots , N$ and $M = 1, 2, \dots$ as follows.
    The              quantities $u_{i}^{( 1 )}$ $( i = 1, \dots , N )$ are defined as
    \begin{equation*}
        u_{i}^{( 1 )} = \begin{cases}
                            0,                               &  i = 1,           \\
                            \tilde{F}_{i} W_{i - 1, i - 1},  &  i = 2, \dots N.
                        \end{cases}
    \end{equation*}
    The quantities $u_{i}^{( M )}$ for $i = 1, \dots , N$ and $M = 2, 3, \dots$ are defined as 
    \begin{gather*}
        u_{i}^{( 2 )} = \sum_{j_{1    } = 1}^{i - 1} W_{i, j_{1}}                                                  W_{j_{1    }, i},                             \\
        u_{i}^{( 3 )} = \sum_{j_{1    } = 1}^{i - 1}
                        \sum_{j_{2    } = 1}^{i - 1} W_{i, j_{1}} W_{j_{1}, j_{2}}                                 W_{j_{2    }, i},                             \\
        u_{i}^{( r )} = \sum_{j_{1    } = 1}^{i - 1}
                        \cdots
                        \sum_{j_{r - 1} = 1}^{i - 1} W_{i, j_{1}} W_{j_{1}, j_{2}} \cdots W_{j_{r - 2}, j_{r - 1}} W_{j_{r - 1}, i},  \qquad  r = 4, 5, \dots .
    \end{gather*}
    For $M = 1, 2, \dots$, the quantities $u _{1}^{( M )}$ and $\tilde{g}_{1}^{( M )}$ coincide since they are zero.
    From the definitions,  the quantities $u _{i}^{( 1 )}$ and $\tilde{g}_{i}^{( 1 )}$ coincide for $i = 2, \dots , N$.
    Then, when we show that a recurrence relation which the quantities $u_{i}^{( M )}$ for $i = 2, \dots , N$ and $M = 2, 3, \dots$ satisfy
                            has the same form as (\ref{rew_gh}) in Remark \ref{Rmk_recur_g}, proof completes.
    We consider the cases where $M = 2$.
    It holds $\tilde{g}_{1}^{( 2 )} = u_{1}^{( 2 )}$ as mentioned above.
    The quantities $u_{i}^{( 2 )}$ for $i = 2, \dots , N$ satisfy
    \begin{align}
        u_{i}^{( 2 )} &=               \sum_{j = 1}^{i - 1}               W_{i    , j} W_{j, i    }
                       =               \sum_{j = 1}^{i - 1} \tilde{F}_{i} W_{i - 1, j} W_{j, i - 1}                                                     \notag         \\
                      &= \tilde{F}_{i} \sum_{j = 1}^{i - 2}               W_{i - 1, j} W_{j, i - 1} + \tilde{F}_{i} W_{i - 1, i - 1} W_{i - 1, i - 1}   \notag         \\
                      &= \tilde{F}_{i} u_{i - 1}^{( 2 )}                                            +               W_{i - 1, i - 1} u_{i}^{( 1 )}      \label{rr_M2}
    \end{align}
    from (\ref{WW_dd}) in Lemma \ref{LemmaEntW} and the definition of $u_{i}^{( 1 )}$.
    Then, the recurrence relations (\ref{rew_gh}) for $M = 2$ in Remark \ref{Rmk_recur_g} and (\ref{rr_M2}) have the same form.
    Thus, $\tilde{g}_{i}^{( 2 )}$ and $u_{i}^{( 2 )}$ for $i = 1, \dots , N$ coincide.
    
    Next, we consider the cases where $i = 2$.
    As mentioned above, it holds $\tilde{g}_{2}^{( M )} = u_{2}^{( M )}$ for $M = 1$ and $2$.
    Then let us consider cases where $M = 3, 4, \dots$.
    It holds $\tilde{g}_{2}^{( M )} = W_{1, 1} \tilde{g}_{2}^{( M - 1 )}$ from (\ref{rew_gh}) in Remark \ref{Rmk_recur_g} since $\tilde{g}_{1}^{( l )} = 0$ $( l = 2, \dots , M )$.
    It holds $u_{2}^{( M )} = W_{2, 1} ( W_{1, 1} )^{M - 2} W_{1, 2} = W_{1, 1} u_{2}^{( M - 1 )}$ from the definition.
    Then, for $M = 1, 2, \dots$, the quantities $\tilde{g}_{2}^{( M )}$ and $u_{2}^{( M )}$ coincide.

    Lastly, let us consider cases where $M$, $N$ and $i$ are $M \geq 3$, $N \geq 3$ and $i = 3, \dots , N$, respectively.
    Let us divide the products $W_{i, j_{1}} \cdots W_{j_{M - 1}, i}$ which appear in the definition of $u_{i}^{( M )}$ into $M$ sets $S_{0}, S_{1}, \dots , S_{M - 1}$.
    The way of division is as follows.
    If     the suffixes   $j_{1}, \dots , j_{M - 1}$ are not $i - 1$, then the product belongs to the set $S_{0}$.
    If     the suffix     $j_{M - 1}$ is $i - 1$, then the product belongs to the set $S_{M - 1}$.
    If     the suffix     $j_{k}$                                is     $j_{k} = i - 1$ for some $k$ $( 1 \leq k \leq M - 2 )$
       and the suffix(es) $j_{l}$ for $l = k + 1, \dots , M - 1$ is/are $j_{l} < i - 1$,
    then the product belongs to the set $S_{k}$.
    Let the summation of the products which belong to the set $S_{k}$ $( k = 0, 1, \dots , M - 1 )$ be denoted by $\alpha _{i}^{( M , k )}$.
    Then, the summation $u_{i}^{( M )}$ is given as
    \begin{equation}  \label{u_alpha}
        u_{i}^{( M )} = \sum_{k = 0}^{M - 1} \alpha _{i}^{( M , k )}.
    \end{equation}
    The summation $\alpha _{i}^{( M, 0 )}$ is written as
    \begin{gather*}
        \alpha _{i}^{( 3, 0 )} = \sum_{j_{1    } = 1}^{i - 2}
                                 \sum_{j_{2    } = 1}^{i - 2} W_{i    , j_{1}} W_{j_{1}, j_{2}}                                 W_{j_{2    }, i},                             \\
        \alpha _{i}^{( r, 0 )} = \sum_{j_{1    } = 1}^{i - 2}
                                 \cdots
                                 \sum_{j_{r - 1} = 1}^{i - 2} W_{i    , j_{1}} W_{j_{1}, j_{2}} \cdots W_{j_{r - 2}, j_{r - 1}} W_{j_{r - 1}, i},  \qquad  r = 4, 5, \dots .
    \end{gather*}
    Since it holds $W_{i, j_{1}} W_{j_{M - 1}, i} = \tilde{F}_{i} W_{i - 1, j_{1}} W_{j_{M - 1}, i - 1}$ from (\ref{WW_dd}) in Lemma \ref{LemmaEntW}, we have
    \begin{equation}  \label{alpha_0}
        \alpha _{i}^{( M, 0 )} = \tilde{F}_{i} u_{i - 1}^{( M )}.
    \end{equation}
    The summation $\alpha _{i}^{( M, M - 1 )}$ is written as
    \begin{gather*}
        \alpha _{i}^{( 3, 2     )} = \sum_{j_{1    } = 1}^{i - 1} W_{i, j_{1}}                         W_{j_{1    }, i - 1} W_{i - 1, i},                             \\
        \alpha _{i}^{( 4, 3     )} = \sum_{j_{1    } = 1}^{i - 1}
                                     \sum_{j_{2    } = 1}^{i - 1} W_{i, j_{1}} W_{j_{1}, j_{2}}        W_{j_{2    }, i - 1} W_{i - 1, i},                             \\
        \alpha _{i}^{( r, r - 1 )} = \sum_{j_{1    } = 1}^{i - 1}
                                     \cdots
                                     \sum_{j_{r - 2} = 1}^{i - 1} W_{i, j_{1}} W_{j_{1}, j_{2}} \cdots W_{j_{r - 2}, i - 1} W_{i - 1, i},  \qquad  r = 5, 6, \dots .
    \end{gather*}
    It holds $W_{j_{M - 2}, i - 1} W_{i - 1, i} = W_{j_{M - 2}, i} W_{i - 1, i - 1}$ from (\ref{WW_ud}) in Lemma \ref{LemmaEntW}.
    Then, we have
    \begin{equation}  \label{alpha_Mm1}
        \alpha _{i}^{( M, M - 1 )} = W_{i - 1, i - 1} u_{i}^{( M - 1 )}.
    \end{equation}
    The summation $\alpha _{i}^{( M, 1 )}$ is written as
    \begin{gather*}
        \alpha _{i}^{( 3, 1 )} = \sum_{j_{2    } = 1}^{i - 2} W_{i, i - 1}                         W_{i - 1    , j_{2    }} W_{j_{2    }, i},                             \\
        \alpha _{i}^{( 4, 1 )} = \sum_{j_{2    } = 1}^{i - 2}
                                 \sum_{j_{3    } = 1}^{i - 2} W_{i, i - 1} W_{i - 1, j_{2}}        W_{j_{2    }, j_{3    }} W_{j_{3    }, i},                             \\
        \alpha _{i}^{( r, 1 )} = \sum_{j_{2    } = 1}^{i - 2}
                                 \cdots
                                 \sum_{j_{r - 1} = 1}^{i - 2} W_{i, i - 1} W_{i - 1, j_{2}} \cdots W_{j_{r - 2}, j_{r - 1}} W_{j_{r - 1}, i},  \qquad  r = 5, 6, \dots .
    \end{gather*}
    Since it holds $                        W_{i    , i - 1} W_{j_{M - 1}, i    }
                    = \tilde{F}_{i}         W_{i - 1, i - 1} W_{j_{M - 1}, i - 1}
                    =        u _{i}^{( 1 )}                  W_{j_{M - 1}, i - 1}$ from (\ref{WW_dd}) in Lemma \ref{LemmaEntW} and the definition of $u_{i}^{( 1 )}$, we have
    \begin{equation}  \label{alpha_1}
        \alpha _{i}^{( M, 1 )} = u_{i}^{( 1 )} u_{i - 1}^{( M - 1 )}.
    \end{equation}
    In the case of $M \geq 4$, the summation $\alpha _{i}^{( M, k )}$ $( k = 2, \dots, M - 2 )$ is given as 
    \begin{equation*}
        \alpha _{i}^{( M, k )} = \sum_{j_{1} = 1}^{i - 1} \cdots \sum_{j_{k - 1} = 1}^{i - 1} \sum_{j_{k + 1} = 1}^{i - 2} \cdots \sum_{j_{M - 1} = 1}^{i - 2}
                                 W_{i, j_{1}} \cdots W_{j_{k - 1}, i - 1} W_{i - 1, j_{k + 1}} \cdots W_{j_{M - 1}, i}.
    \end{equation*}
    Since it holds $W_{j_{k - 1}, i - 1} W_{j_{M - 1}, i} = W_{j_{k - 1}, i} W_{j_{M - 1}, i - 1}$ from (\ref{WW_ud}) in Lemma \ref{LemmaEntW}, we have
    \begin{equation}  \label{alpha_Mid}
        \alpha _{i}^{( M, k )} = u_{i}^{( k )} u_{i - 1}^{( M - k )}.
    \end{equation}
    From (\ref{u_alpha}), (\ref{alpha_0}), (\ref{alpha_Mm1}), (\ref{alpha_1}) and (\ref{alpha_Mid}), we obtain
    \begin{align}
        u_{i}^{( M )} &=                                                                        \sum_{k = 0}^{M - 1} \alpha _{i}^{( M , k )}
                       = \tilde{F}_{i} u_{i - 1}^{( M )} + W_{i - 1, i - 1} u_{i}^{( M - 1 )} + \sum_{k = 1}^{M - 2} u_{i    }^{( k )} u_{i - 1}^{( M - k )}   \notag            \\
                      &= \tilde{F}_{i} u_{i - 1}^{( M )} + W_{i - 1, i - 1} u_{i}^{( M - 1 )} + \sum_{k = 2}^{M - 1} u_{i - 1}^{( k )} u_{i    }^{( M - k )}.  \label{rr_uhigh}
    \end{align}
    Then, the recurrence relations (\ref{rew_gh}) in Remark \ref{Rmk_recur_g} and (\ref{rr_uhigh}) have the same form.
    Thus, $\tilde{g}_{i}^{( M )}$ and $u_{i}^{( M )}$ for $i = 3, \dots , N$ and $M = 3, 4, \dots$ coincide. ~ $\square$
    
\subsection*{Proof of Theorem \ref{Th_trans_gh}}

    Let us introduce quantities $\zeta _{i}^{( M )}$ for $i = 1, \dots , N$ and $M = 1, 2, \dots$ defined as
    \begin{align}
        &\zeta _{1}^{( 1 )} = 0,                                      & &                   & &                   \label{zeta11}    \\
        &\zeta _{i}^{( 1 )} = \tilde{F}_{i}       h_{i - 1}^{( 1 )},  & &i = 2, \dots , N,  & &                   \label{zetai1}    \\
        &\zeta _{i}^{( M )} = \frac{1}{M !} \cdot h_{i    }^{( M )},  & &i = 1, \dots , N,  & &M = 2, 3, \dots .  \label{zetaHigh}  
    \end{align}
    When we show that a recurrence relation which the quantities $\zeta _{i}^{( M )}$ satisfy
                      coincides with that of the quantities $\tilde{g}_{i}^{( M )}$ in Definition \ref{def_tg}, proof completes.

    It is obvious that $\zeta _{1}^{( M )} = \tilde{g}_{1}^{( M )}$ for $M = 1, 2, \dots$ since they are zero from their definitions.
    
    Hereafter, let $i$ be $i = 2, \dots , N$ in this proof.

    The quantities $h_{i - 1}^{( 1 )}$ coincide with the diagonal entry $W_{i - 1, i - 1}$ of $W$ from Remark \ref{Rmk_diff_w}.
    Then, it holds
    \begin{equation*}
        \zeta _{i}^{( 1 )} = \tilde{g}_{i}^{( 1 )}
    \end{equation*}
    from the definition of $\tilde{g}_{i}^{( 1 )}$.

    Let us consider the case where $M = 2$.
    From (\ref{zetaHigh}), the recurrence relation for $h_{i}^{( 2 )}$ in Theorem \ref{Th_G} is rewritten as
    \begin{equation*}
        2 \zeta _{i}^{( 2 )} = \tilde{F}_{i} \left( 2 \zeta _{i - 1}^{( 2 )} + 2 \left( h_{i - 1}^{( 1 )} \right) ^{2} \right) .
    \end{equation*}
    From (\ref{zetai1}) and the coincidence $h_{i - 1}^{( 1 )} = W_{i - 1, i - 1}$, we have
    \begin{equation*}
        \zeta _{i}^{( 2 )} = \tilde{F}_{i} \zeta _{i - 1}^{( 2 )} + W_{i - 1, i - 1} \zeta _{i}^{( 1 )}.
    \end{equation*}
    This recurrence relation has the same form with the recurrence relation (\ref{rew_gh}) for $M = 2$ in Remark \ref{Rmk_recur_g}.

    Let us consider the case where $M = 3, 4, \dots$.
    The recurrence relation (\ref{hih}) for $h_{i}^{( M )}$ in Theorem \ref{Th_G} is rewritten as
    \begin{equation*}
        h_{i}^{( M )} =     \tilde{F}_{i}                                     h_{i - 1}^{( M     )}
                        + M \tilde{F}_{i}                   h_{i - 1}^{( 1 )} h_{i - 1}^{( M - 1 )}
                        + M                                 h_{i - 1}^{( 1 )} h_{i    }^{( M - 1 )}
                        + \sum_{k = 2}^{M - 2} {}_{M} C_{k} h_{i - 1}^{( k )} h_{i    }^{( M - k )}.
    \end{equation*}
    On the second and the third terms in the right-hand-side, it holds
    \begin{gather*}
            M \tilde{F}_{i} h_{i - 1}^{( 1 )} h_{i - 1}^{( M - 1 )}
          + M               h_{i - 1}^{( 1 )} h_{i    }^{( M - 1 )}                                                                  \\
        =   M                                   \zeta _{i}^{( 1     )} \cdot ( M - 1 ) ! \cdot \zeta _{i - 1}^{( M - 1 )}
          + M                  W_{i - 1, i - 1}                        \cdot ( M - 1 ) ! \cdot \zeta _{i    }^{( M - 1 )}            \\
        =   M ! \cdot \left(   W_{i - 1, i - 1}                                                \zeta _{i    }^{( M - 1 )}
                             +                  \zeta _{i}^{( 1 )}                             \zeta _{i - 1}^{( M - 1 )} \right) .
    \end{gather*}
    Then, we have
    \begin{align*}
        M ! \cdot \zeta _{i}^{( M )} = &  \tilde{F}_{i} \cdot M ! \cdot                                                                     \zeta _{i - 1}^{( M     )}
                                        +                     M ! \cdot \left( W_{i - 1, i - 1} \zeta _{i}^{( M - 1 )} + \zeta _{i}^{( 1 )} \zeta _{i - 1}^{( M - 1 )} \right)  \\
                                       &+ \sum_{k = 2}^{M - 2} ~_M C_{k} \cdot k ! \cdot \zeta _{i - 1}^{( k )} \cdot ( M - k ) ! \cdot \zeta _{i}^{( M - k )}
    \end{align*}
    from the definition (\ref{zetaHigh}).
    Since it holds $~_M C_{k} \cdot k ! \cdot ( M - k ) ! = M !$, we immediately have
    \begin{equation*}
        \zeta _{i}^{( M )} =                        \tilde{F}_{i}                  \zeta _{i - 1}^{( M )}
                             +                                    W_{i - 1, i - 1}                        \zeta _{i}^{( M - 1 )}
                             + \sum_{k = 2}^{M - 1}                                \zeta _{i - 1}^{( k )} \zeta _{i}^{( M - k )}.
    \end{equation*}
    This recurrence relation has the same form as (\ref{rew_gh}) for $M = 3, 4, \dots$ in Remark \ref{Rmk_recur_g}. ~ $\square$

\subsection*{Proof of Theorem \ref{Th_NewFormula}}
    
    From Remarks \ref{Rmk_diff_w} and \ref{Rmk_diag_NewFormula}, it holds
    \begin{equation}  \label{tGHh_diag}
        \tilde{G}_{i}^{( 1 )} = H_{i}^{( 1 )} = h_{i}^{( 1 )},  \qquad  i = 1, \dots , N.
    \end{equation}
    Substituting $h_{i}^{( 1 )} = \tilde{G}_{i}^{( 1 )}$ $( i = 1, \dots , N )$ into (\ref{h11}) and (\ref{hi1}) in Theorem \ref{Th_G}
    and defining $\tilde{g}_{i}^{( 1 )}$ $( i = 1, \dots , N )$ as (\ref{tg11}) and (\ref{tgi1}),
    we readily have (\ref{tG11}) and (\ref{tGi1}).
    
    Next, we substitute (\ref{trans_gh}) in Theorem \ref{Th_trans_gh} into (\ref{h1h}) and (\ref{hih}) in Theorem \ref{Th_G}.
    We readily have (\ref{tg1h}).
    And we have
    \begin{align*}
                                                                                                M       ! \cdot \tilde{g}_{i    }^{( M     )}
        =&  \tilde{F}_{i} \left(                                                                M       ! \cdot \tilde{g}_{i - 1}^{( M     )}
                                 +                  M   \cdot        h _{i - 1}^{( 1 )} \cdot ( M - 1 ) ! \cdot \tilde{g}_{i - 1}^{( M - 1 )} \right)  \\
         &+                                         M   \cdot        h _{i - 1}^{( 1 )} \cdot ( M - 1 ) ! \cdot \tilde{g}_{i    }^{( M - 1 )}          \\
         &+ \sum_{k = 2}^{M - 2} {}_{M} C_{k} \cdot k ! \cdot \tilde{g}_{i - 1}^{( k )} \cdot ( M - k ) ! \cdot \tilde{g}_{i    }^{( M - k )}
    \end{align*}
    for $i = 2, \dots , N$ and $M = 2, 3, \dots$.
    Then, it holds
    \begin{equation*}
                                                           \tilde{g}_{i    }^{( M     )}
        =   \tilde{F}_{i}                                  \tilde{g}_{i - 1}^{( M     )}
          + \tilde{F}_{i}               h _{i - 1}^{( 1 )} \tilde{g}_{i - 1}^{( M - 1 )}
          +                             h _{i - 1}^{( 1 )} \tilde{g}_{i    }^{( M - 1 )}
          + \sum_{k = 2}^{M - 2} \tilde{g}_{i - 1}^{( k )} \tilde{g}_{i    }^{( M - k )}.
    \end{equation*}
    The second term of the right-hand-side is equal to $\tilde{g}_{i}^{( 1 )} \tilde{g}_{i - 1}^{( M - 1 )}$
    since it holds $\tilde{F}_{i} h_{i - 1}^{( 1 )} = \tilde{F}_{i} \tilde{G}_{i - 1}^{( 1 )} = \tilde{g}_{i}^{( 1 )}$ from (\ref{tgi1}) and (\ref{tGHh_diag}).
    From (\ref{tGHh_diag}), the third term of the right-hand-side is equal to $\tilde{G}_{i - 1}^{( 1 )} \tilde{g}_{i}^{( M - 1 )}$.
    Then, we have (\ref{tgih}).

    Next, we substitute (\ref{trans_gh}) in Theorem \ref{Th_trans_gh} and (\ref{Transf_HtG}) into (\ref{Hhigh}) in Theorem \ref{Th_G}.
    We have
    \begin{align*}
                                                                                              ( M - 1     ) ! \cdot \tilde{G}_{i}^{( M     )}
        =&                                                                                      M           ! \cdot \tilde{g}_{i}^{( M     )}
          +                        ( M - 1 )                             h _{i}^{( 1 )} \cdot ( M - 2     ) ! \cdot \tilde{G}_{i}^{( M - 1 )}   \\
         &+ \sum_{k = 2}^{M - 1} {}_{M - 1} C_{k} \cdot k ! \cdot \tilde{g}_{i}^{( k )} \cdot ( M - 1 - k ) ! \cdot \tilde{G}_{i}^{( M - k )}
    \end{align*}
    for $i = 1, \dots , N$ and $M = 2, 3, \dots$.
    Then, it holds
    \begin{equation*}
                                                       \tilde{G}_{i}^{( M     )}
        =                 M      \tilde{g}_{i}^{( M )}
          +                             h _{i}^{( 1 )} \tilde{G}_{i}^{( M - 1 )}
          + \sum_{k = 2}^{M - 1} \tilde{g}_{i}^{( k )} \tilde{G}_{i}^{( M - k )}
    \end{equation*}
    for $i = 1, \dots , N$ and $M = 2, 3, \dots$.
    Since the second term of the right-hand-side is equal to $\tilde{G}_{i}^{( 1 )} \tilde{G}_{i}^{( M - 1 )}$ from (\ref{tGHh_diag}), we obtain (\ref{tGhigh}).
    
    Lastly, substituting (\ref{Transf_HtG}) into (\ref{TraceByH}) in Theorem \ref{Th_G}, we readily have (\ref{TraceByTG}). ~ $\square$

\end{document}